\documentclass[10pt]{article}
\usepackage{amsfonts}                   
\usepackage{amssymb}
\usepackage{amsmath}

\usepackage{verbatim}

\def\Prob{\mathbb{P}}
\def\Pr{\mathbb{Q}}
\def\Ev{\mathbb{E}}
\def\stand{W_{ST}^d}
\def\double{W_{DD}^d}
\def\hit{W_H^d}

\def\deck{\mbox{\it deck}}
\begin{document}
\title{Finding Blackjack's Optimal Strategy in Real-time and Player's Expected Win}
\author{Jarek Solowiej}
\date{\today}
\maketitle

\footnotetext[1]{jarek@math.albany.edu}
\begin{abstract}
We describe the probability theory behind a casino game, blackjack, and 
the procedure to compute the optimal strategy for a deck of arbitrary 
cards and player's expected win given that he follows the optimal strategy. 
The exact blackjack probabilities are used, in contrast 
to approximate probabilities used by Baldwin et al.~\cite{BCMM} 
or Monte Carlo methods.
We distinguish between two probability measures $\Prob$ and $\Pr$;
$\Prob$ is used to compute dealer's probabilities and
$\Pr$ is used to compute player's expectations.
The implementation is described in pseudo-C++. 
The program is fast enough to deal with any blackjack's 
hand in a matter of seconds. 
\end{abstract}

\section*{ }
The main rules of blackjack are identical for all casinos, but 
they differ on details. And it would be too cumbersome to deal 
with all the variation of those. 
So this paper provides only the general framework, 
which, by the way, is easily adaptable to any particular set of blackjack's rules.
The reader can see \cite{BCMM} for explanation of the rules. 

We distinguish between hard and soft blackjack's
hands; a hand with an ace counted as 11 is soft, otherwise we call it hard.
For example, ace and five is a soft sixteen  ($1+5=16s$); 
ten and six is a hard sixteen ($10+6=16h$); 
or ace and ace is a soft twelve ($1+1=12s$).

\section*{Dealer's Probabilities}
During the first stage of a game, 
dealer gets two cards, one face up and one face down. 
We can assume that 
his face-down card is still in a deck.
Let $D_d$ denote dealer's total given that his face-up card is $d$ 
(in this paper $d$ always denotes dealer's face-up card), where $d=1$ means
an ace (it can be counted as 1 or 11), 
$d=2$ means a deuce, $\ldots$, 
and $d=10$ means a ten or any face card (we don't distinguish between a ten and any face card and
we refer to them as a ten). 

The blackjack's rules force dealer to hit
anything below 17, so dealer's total can be 
17, 18, 19, 20, 21 (any combination of cards adding to 21 
but ten and ace), 
a natural (ten and ace), 
and a bust (anything above 21). To simplify the
notation we denote a natural by 22 and a bust by 23.

At this point let's assume that we know
dealer's probabilities  (see the appendix for the algorithm),
$$
\Prob[D_d=k]
\quad\mbox{for}\quad k=17,\dots,22,\quad\mbox{and}\quad d=1,2,3,\ldots,10,
$$
and $\Prob[D_d=23]=1-\sum_{k=17}^{22}\Prob[D_d=k]$. 
We assume that there are enough cards in a deck
so $\Prob[D_d=k]=0$ for $k<17$,
and therefore the distribution,  $\Prob[D_d=17],\ldots,\Prob[D_d=23]$,
is not degenerate.

Since dealer checks for a natural (if $d=1$ or $d=10$) before player
makes any decision we exclude this possibility by 
conditioning on the event $D_d\neq 22$. This
provides us with a new probability 
measure $\Pr$ such that $\Pr[D_d=22]=0$.
By definition, given $d=1,2,\ldots,10$, 
$$
\Pr[D_d=k]=\Prob[D_d=k|D_d\neq 22]=\frac{\Prob[D_d=k]}{\Prob[D_d\neq 22]}
\quad\mbox{for}\quad k=17,\dots,21,
$$
where the trivial case  $\Prob[D_d\neq 22]= 0$ 
or $\Prob[D_d=22]=1$ is excluded.
Then $\Pr[D_d=22]=0$ and we set
$\Pr[D_d=23]=1-\sum_{k=17}^{21}\Pr[D_d=k]$.
Note that $\Prob=\Pr$ for $d\neq 1, 10$, since 
$\Prob[D_d\neq 22]=1$ for $d\neq 1, 10$.

The probability measure $\Pr$ 
and player's cards determine the optimal strategy. For example,
given that player's total equals 19, the player is interested
in $\Pr[D_d=20]+\Pr[D_d=21]$, since this is the probability that
he loses if he stands. 
\begin{table}[ht]
\caption{Dealer's probabilities
using $\Pr$ for one deck (dealer stands on soft 17; the numbers are cut after
five digits).}
\begin{center}
\begin{tabular}{|c|c|c|c|c|c|c|}
\hline 
$\Pr$ & 17    &    18   &     19  &    20   &    21   &    bust  \\ \hline
 2  &0.13897&  0.13176&  0.13181&  0.12394&  0.12052&  0.35297 \\
 3  &0.13030&  0.13094&  0.12376&  0.12334&  0.11604&  0.37559 \\
 4  &0.13097&  0.11416&  0.12067&  0.11628&  0.11509&  0.40280 \\
 5  &0.11968&  0.12348&  0.11690&  0.10469&  0.10632&  0.42890 \\
 6  &0.16694&  0.10645&  0.10719&  0.10070&  0.09787&  0.42082  \\
 7  &0.37234&  0.13858&  0.07733&  0.07889&  0.07298&  0.25985 \\
 8  &0.13085&  0.36298&  0.12944&  0.06828&  0.06979&  0.23862 \\
 9  &0.12188&  0.10392&  0.35739&  0.12225&  0.06110&  0.23344 \\
10  &0.12415&  0.12248&  0.12442&  0.35686&  0.03956&  0.23249 \\
 1  &0.18378&  0.19089&  0.18868&  0.19169&  0.07513&  0.16981 \\
\hline
\end{tabular}
\end{center}
\end{table}

\section*{Deck and Cards' Probabilities}
Now, we want to know how the fact that dealer doesn't have
a natural affects cards' probabilities.

At any time the content of a deck is described by ten numbers $(a_1,a_2,\ldots,a_{10})$,
where $a_1$ denotes the number of aces, 
$a_2$ the number of deuces, $a_3$ the number of three's, 
$\ldots$, and $a_{10}$ 
the number of ten's. 
We refer to this deck as $\deck$.
When  a card with a face value $k=1,2,\ldots,10$ is drawn from $\deck$, then 
$$
\Prob[k]=a_k/t,
\quad\mbox{where}\quad t=a_1+a_2+\cdots+a_{10}.
$$
Now suppose that dealer's face-up card is a ten or an ace, 
but he doesn't have a natural. 
Despite the fact that we use the same deck,
the probabilities are different.
If $d=10$, then
$$
\Pr[1]=\frac{a_1}{t-1},\quad \Pr[k]=\frac{a_k}{t-1}\frac{t-a_1-1}{t-a_1}
\quad\mbox{for}\quad k=2,\ldots,10,
$$
or if $d=1$, then
$$
\Pr[10]=\frac{a_{10}}{t-1},\quad \Pr[k]=\frac{a_k}{t-1}\frac{t-a_{10}-1}{t-a_{10}}
\quad\mbox{for}\quad k=1,\ldots,9.
$$
And as we noted before $\Prob[k]=\Pr[k]$ for $k=1,2,\ldots,10$ if $d\neq 1,10$.
See the appendix for derivation of these probabilities.

\section*{Player's Options}
In this section we assume that dealer doesn't have a natural. A player
with two cards has three or four options (this depends on casino's rules).
We use $W_{ST}^d$ to denote  player's win if he stands and dealer's
face-up card is $d$; $W_{DD}^d$, he double downs;  $W_H^d$, 
he hits;  $W_S^d$, he splits.
Given player's hand $p$, 
the optimal strategy is determined by looking at expected win 
for each of the available options:
$\Ev[W_{ST}^d|p]$, $\Ev[W_{DD}^d|p]$, $\Ev[W_H^d|p]$, or $\Ev[W_S^d|p]$.
An option with the highest expected win is optimal.

If a player has a natural, 
then he wins one and a half of his original bet, and the game is over so
we exclude this possibility in the computations below.

\subsubsection*{Stand}
Given a player with total $p\le 21$ (no naturals);
if he stands, then the three possible outcomes are as
follows 
\begin{itemize}
\item if $D_d<p$ or $D_d=23$, then he wins one dollar;
\item if $D_d>p$ and $D_d\neq 23$, then he loses one dollar;
\item if $D_d=p$, we have a tie, nothing happens. 
\end{itemize}
So player's expected win equals
$$
\Ev[\stand|p,\deck]=\Pr[D_d<p]+\Pr[D_d=23] -\Pr[D_d>p, D_d\neq 23],
$$
where $\Ev[\ \cdot\ |p,\deck]$ denotes the expected win  
under the condition 
that his total equals $p$ and
the cards are coming from  $\deck$. 
To simplify the notation we usually skip the conditioning on $\deck$, moreover
when dealing with expressions of the type $\Pr[i]\Ev[\ \cdot\ |p+i]$, it is understood
that the expected value is conditioned on $\deck-\{i\}$.

\subsubsection*{Double Down}
Given a player with two cards,
doubling down means that
he doubles his bet and gets only one additional card. 
By considering all the possibilities for this extra card, we get
$$
\Ev[\double|p]=2\sum_{i=1}^{10}\Pr[i]\Ev[\stand|p+i].
$$

\subsubsection*{Hit}
Given a player with total $p$,
when hitting he gets one card and has an option to stand or take 
another hit; we can assume that the player doesn't hit 21, since
$\Ev[\hit|21]\le \Ev[\stand|21]$ in any situation. Therefore
$$
\Ev[\hit|p]=\sum_{i=1}^{10}\Pr[i]
\begin{cases}
\max(\Ev[\stand|p+i],\Ev[\hit|p+i]) &\mbox{if}\quad p+i<21; \\
\Ev[\stand|p+i]                            &\mbox{otherwise.}
\end{cases}
$$
It is a finite recursion, since   the longest sequence of player's cards has
length 20.\footnote{1+1+1+1+1+1+1+1+1+1+2+1+1+1+1+1+1+1+1+1=21.} 
We can introduce an extra variable that controls 
the recursion level; 13 is needed to compute $\Ev[\hit|1+1]$, 
but other hands need lower level of recursion, given desired  accuracy, see the implementation.

\subsubsection*{Split}
If a player has two cards with the same  value, he can split them into
two hands (adding extra bet). Casinos' rules  
on player's option after a split vary; they differ 
on  number of splits allowed for non-aces and aces, 
and the ability to double down.  
Despite that,  the following is true for any split
$$
\Ev[W_S^d|p+p]=2\Ev[W_{\mbox{play}}^d|p],
$$
where ``play'' means all the options available to a player
after splitting a hand. This formula is true because
after split there are two identical hands (they are not independent but
this doesn't concern us since expected value of 
sum equals sum of expected values).

When splitting aces each hand gets only one additional card (blackjack's rule), so 
$$
\Ev[W_{S1}^d|1+1]=2\sum_{i=1}^{10}\Pr[i]\Ev[\stand|1+i],
$$
where $S1$ means only one split. 
If a casino allows for re-splitting aces then 
``play'' means only one additional card but two aces
can be split again, so
$$
\Ev[W_S^d|1+1]
=2\Pr[1]\max(\Ev[\stand|1+1],\Ev[W_{S1}^d|1+1])+2\sum_{i=2}^{10}\Pr[i]\Ev[\stand|1+i].
$$
See the implementation for the general case.

\begin{table}[ht]
\caption{The expectations for two decks 
(dealer stands on soft 17, double down after split is allowed, 
re-splitting aces and non-ace pairs is allowed; the numbers are cut after six digits).}
\begin{center}
\begin{tabular}{|c|r|r|r|r|r|c|}
\hline
$p$ & $\Ev[W_{ST}^6|p]$ & $\Ev[W_H^6|p]$ & $\Ev[W_{DD}^6|p]$ & $\Ev[W_S^6|p]$ & Action \\ \hline
$2+10$   & $ -0.156818$ &    $   -0.165123$ & $     -0.330246$  & & stand \\
$3+10$   & $ -0.155641$ &    $  -0.232503 $ & $      -0.465006$ & & stand\\
$4+10$   & $ -0.154544$ &    $  -0.304424 $ & $    -0.608848 $ &  & stand\\
$5+10$   & $ -0.153729$ &    $  -0.376364 $ & $    -0.752728 $ &  & stand\\
$6+10$   & $ -0.165609$ &    $  -0.414113 $ & $    -0.828226 $ & & stand \\
$7+10$   & $ 0.001024$ &     $ -0.496273 $  & $   -0.992546 $ &  &stand\\
$8+10$   & $  0.276027$ &    $   -0.597068$ &  $    -1.194140$ & &stand\\
$9+10$   & $  0.490271$  &   $   -0.714945$  & $    -1.429890$ & &stand\\
$10+ 10$  & $ 0.700605$    &  $  -0.849453 $  & $   -1.698910$    & $0.569494$&stand\\
$1+ 1 $   &  $-0.129268 $  &  $  0.192311 $   &   $ 0.213109 $    & $  0.836235$ &split\\
$1+ 2$    & $ -0.134355  $ & $   0.164810 $   &    $     0.204564$ & &double \\
$1+ 3 $   & $ -0.133179 $  &  $  0.142659 $   &  $  0.200079 $& &double\\
$1+ 4 $   & $ -0.132096 $  &  $  0.118918 $   &  $  0.189631$ & &double\\
$1+ 5$    & $ -0.131183 $  &  $  0.107088 $   &  $  0.197579$ & &double\\
$1+ 6$    & $ 0.012003 $  &   $ 0.131284 $   &  $  0.262569$ &  &double\\
$1+ 7$    & $ 0.273910 $   &   $   0.192289 $ & $     0.384579 $& &double\\
$1+ 8$    & $ 0.489571 $   &   $ 0.240709 $   & $   0.481418 $ & &stand\\
$1+ 9$    & $ 0.699584$    &   $ 0.284227 $   & $   0.568454$ & &stand\\
$1+ 10$   &  $1.500000$    &$ 0.337395$       & $0.674791$  & & stand\\ \hline
\end{tabular}
\end{center}
\end{table}

\subsection*{Optimal Strategy and It's Payoff}
Given player's two-card hand with total $p$, the optimal strategy is given
by the action with the highest expectation. For example, if
$$
\Ev[\double|p]= \max(\Ev[\stand|p],\Ev[\double|p],\Ev[\hit|p],\Ev[W_S^d|p]),
$$
then a player should double down. Or for player's hand 
with at least three cards, we compare
$\Ev[\stand|p]$ to $\Ev[\hit|p]$.

Now, we focus on an arbitrary hand,
assuming that player follows the optimal strategy. We want to
compute his average win $W$.
Let $W_d$ denote player's win given that dealer's face-up card is $d$, then
$\Ev[W]=\sum_{d=1}^{10}\Prob[d]\Ev[W_d]$.
For $d\neq 1,10$, we have
$$
\Ev[W_d]=\sum_{i,j=1}^{10}\Pr[i,j]\Ev[W_d|i+j],
\quad\mbox{where}\quad \Ev[W_d|i+j]=1.5\mbox{ if } i+j=21, 
$$
and
$\Ev[W_d|i+j]=\max(\Ev[\stand|i+j],\Ev[\hit|i+j],\Ev[\double|i+j],\Ev[W_S^d|i+j])$
if $i+j<21$.

There are two special cases,  $\Ev[W_1]$ and $\Ev[W_{10}]$, where we need to consider
the possibility that dealer has a natural.
In order to compute $\Ev[W_1]$ we consider two cases: dealer has a natural and dealer 
does not have a natural,
\begin{align*}
\Ev[W_1]
&=\Prob[D_1=22]\Ev[W_1|D_1=22]+\Prob[D_1\neq 22]\Ev[W_1|D_1\neq 22] \\
&=\Prob[10]\Ev[W_1|D_1=22]+(1-\Prob[10])\Ev[W_1|D_1\neq 22] \\
&=-\Prob[10](1-\Prob[1,10|\deck-\{1,10\}])+(1-\Prob[10])\sum_{i,j=1}^{10}\Pr[i,j]\Ev[W_1|i+j]. 
\end{align*}
Note that $\Prob[1,10]$ means the probability of a natural.
Moreover, we compute $\Ev[W_{10}]$ using the same method.

\begin{table}[ht]
\caption{The expected wins, $\Ev[W]$, as a percentage of the original bet for different
set of rules and different number of decks (the first column).
The first row specifies the options used: double down after split, re-split aces, re-split non-ace pairs.
For example, $010$ means double down after split is not allowed, re-splitting aces is allowed, and
re-splitting non-ace pairs is not allowed. The card probabilities are constant for the infinite deck.}
\begin{center}
{\small
\begin{tabular}{|c|c|c|c|c|c|c|c|c|}
\hline 
$\Ev[W]$ & 000      & 001       & 010      &  011       &  100       &  101       &   110      &   111     \\ \hline
 1     & $-.6747$ & $-.6585$  & $-.6508$ &  $-.6345$  &  $-.5452$  &  $-.5174$  &  $-.5212$  & $-.4934$  \\
 2     & $-.6876$ & $-.6628$  & $-.6438$ &  $-.6190$  &  $-.5619$  &  $-.5215$  &  $-.5180$  & $-.4776$  \\
 3     & $-.6884$ & $-.6607$  & $-.6369$ &  $-.6093$  &  $-.5645$  &  $-.5202$  &  $-.5131$  & $-.4688$  \\
 4     & $-.6889$ & $-.6598$  & $-.6335$ &  $-.6045$  &  $-.5659$  &  $-.5197$  &  $-.5105$  & $-.4644$  \\
 5     & $-.6893$ & $-.6594$  & $-.6315$ &  $-.6017$  &  $-.5669$  &  $-.5196$  &  $-.5091$  & $-.4618$  \\
 6     & $-.6897$ & $-.6593$  & $-.6303$ &  $-.5999$  &  $-.5677$  &  $-.5197$  &  $-.5083$  & $-.4603$  \\
 7     & $-.6900$ & $-.6592$  & $-.6294$ &  $-.5986$  &  $-.5683$  &  $-.5197$  &  $-.5077$  & $-.4592$  \\
 8     & $-.6902$ & $-.6591$  & $-.6288$ &  $-.5977$  &  $-.5687$  &  $-.5198$  &  $-.5073$  & $-.4583$  \\
       &          &           &          &            &            &            &            &           \\
 $\infty$    & $-.6901$ & $-.6569$  & $-.6223$ &  $-.5891$  &  $-.5702$  &  $-.5186$  &  $-.5025$  & $-.4509$  \\
\hline
\end{tabular}
}
\end{center}
\label{ex-win}
\end{table}

The first thing to notice about Table~\ref{ex-win} 
is that all its numbers are negative. 
So playing a blackjack is a losing proposition (if you don't vary your bets).
But we can do better.
In theory, given a $\deck$, 
if $\Ev[W|\deck]>0$ we bet one dollar; if  $\Ev[W|\deck]\le 0$ we bet zero dollars. 
This strategy has a positive
expected win, although it has to be modified, 
since computing $\Ev[W|\deck]$ at a blackjack's table would be difficult. 
The solution: we estimate $\Ev[W|\deck]$ by paying attention to cards removed from a deck.

First, we compute the effect of card removal on expected win, namely, we compute
the change in expected value if one card is removed,
$$
r_i=\Ev[W|\deck-\{i\}]- \Ev[W|\deck]\quad\mbox{for}\quad i=1,2,\ldots,10,
$$
see Table~\ref{ex-removal}.
It implies  that the cards 1, 8, 9, and 10 (inside the deck) increase player's  expected win, since 
$r_i<0$ implies $\Ev[W|\deck]> \Ev[W|\deck-\{i\}]$;
the cards 2,3,4,5,6, and 7
decrease player's expected win.
So by counting cards from these two groups and weighting them according to Table~\ref{ex-removal}, 
we can estimate $\Ev[W]$ 
and bet accordingly, for more details see Gottlieb~\cite{GG} or Thorp~\cite{T}.

\begin{table}
\caption{Effect of card removal on expected win in percentages,
where the first row gives the number of decks and the first column represents cards to be removed;
double down after split is allowed, re-splitting aces and non-ace pairs is  allowed.}
\begin{center}
{\small
\begin{tabular}{|c|r|r|r|r|r|r|r|r|}
\hline 
$r_i$  & 1       & 2        & 3       &  4           &  5        &  6          &   7       & 8         \\ \hline
 1     & $-.656$ & $-.316$  & $-.208$ &  $-.155$     &  $-.123$  &  $-.102$    &  $-.087$  & $-.076$   \\
 2     & $.363$  & $.187$   & $.126$  &  $.094$      &  $.076$   &  $.064$     &  $.054$   & $.048$    \\
 3     & $.427$  & $.216$   & $.144$  &  $.109$      &  $.088$   &  $.073$     &  $.063$   & $.055$    \\
 4     & $.564$  & $.289$   & $.194$  &  $.147$      &  $.118$   &  $.099$     &  $.085$   & $.075$    \\
 5     & $.733$  & $.373$   & $.250$  &  $.188$      &  $.151$   &  $.126$     &  $.108$   & $.095$    \\
 6     & $.415$  & $.209$   & $.140$  &  $.106$      &  $.085$   &  $.071$     &  $.062$   & $.054$    \\
 7     & $.267$  & $.134$   & $.089$  &  $.066$      &  $.053$   &  $.044$     &  $.038$   & $.033$    \\
 8     & $-.013$ & $-.010$  & $-.009$ &  $-.008$     &  $-.007$  &  $-.005$    &  $-.004$  & $-.004$   \\
 9     & $-.181$ & $-.097$  & $-.068$ &  $-.051$     &  $-.041$  &  $-.034$    &  $-.029$  & $-.025$   \\
 10    & $-.444$ & $-.236$  & $-.162$ &  $-.122$     &  $-.098$  &  $-.082$    &  $-.070$  & $-.061$      \\
\hline
\end{tabular}
}
\end{center}
\label{ex-removal}
\end{table}

\vskip5pt
\noindent{\bf Example.}
Let $\deck$ be composed of two decks with 111 as our blackjack's rules, see Table~\ref{ex-win}.
So we know that $\Ev[W|\deck]=-.004776$. Now, suppose
that the cards $R=\{1,1,2,2,3,4,4,4,5,5,5,7,7,8,10\}$ are removed from the deck. 
The direct computation shows that $\Ev[W|\deck-R]=0.018249$,
whereas using Table~\ref{ex-removal} and linear interpolation (an interpolation using an 
exponential function would be better)  when removing more than
one card of the same value, we estimate 
$$
\Ev[W|\deck-R]\approx 0.01965.
$$

\section*{How fast are these procedures when implemented?}
The execution times of the program implementing the procedures described in this paper 
depends on dealer's face-up card. The table below gives times (in seconds)
needed to compute  $\Ev[\stand|p]$ one thousand times for different dealer's face-up cards
(player's total is negligible, four decks were used). It was done 
on a computer with two 1.2 MHz processors (with 256 kB cache each), Although,
only one processor was used to run the program.

\begin{center}
\begin{tabular}{|c|c|c|c|c|c|c|c|c|c|c|}
\hline 
$d$  &    2 &   3  &   4 &  5   &  6  &    7 &     8&     9&   10  &   1   \\ \hline
 time & 2.24 & 1.36 & 0.82 &  0.49& 0.29&  0.18&  0.11&  0.07&  0.05 &  0.73 \\
\hline
\end{tabular}
\end{center}

Moreover, $\Ev[\double|p]$ computes $\Ev[\stand|\,\cdot\,]$ ten times;
computational complexity of $\Ev[\hit|p]$ depends on recursion level, for example, 
it takes 0.51 seconds to compute $\Ev[W_H^2|10]$ with 9 recursion levels and
0.05 seconds with 2  recursion levels (the outputs are identical 
up to 15 digits after a decimal point).
The most complicated case $\Ev[W_S^2|2+2]$ 
took 173.76 seconds to compute and 1.10 seconds for $\Ev[W_S^2|10+10]$.
A faster  processor with a bigger cache 
should compute any split in a matter of seconds.
Moreover, it takes about 30 minutes to compute $\Ev[W]$, although
some trivial optimization techniques could cut that in half.

\section*{Appendix A: conditional probabilities}
Given $d=10$ and the fact that the dealer doesn't have a natural
we want to compute $\Pr[k]$ for $k=1,2,\ldots,10$. By definition
\begin{align*}
\Pr[k]&=\Prob[k|\mbox{dealer's face-down card $\neq$ ace}] \\
&=\frac{\Prob[k \mbox{ and dealer's face-down card $\neq$ ace}]}
{\Prob[\mbox{dealer's face-down card $\neq$ ace}]}.
\end{align*}
We compute $\Pr[1]$; first the probability in the denominator,
\begin{align*}
\Prob[\mbox{dealer's face-down card $\neq$ ace}]
&=1-\Prob[\mbox{dealer's face-down card $=$ ace}]\\ &=1-a_1/t.
\end{align*}
To compute the probability in the numerator, 
consider all the permutations of the unknown cards 
(dealer's face-down card is considered to be in a deck); 
there are $\frac{t!}{a_1!a_2!\cdots a_{10}!}$
of them. Next consider all the permutation with an ace as 
the first and a non-ace as the last card (dealer's face-down card). 
We count them  considering nine choices for the last card,
\begin{align*}
&\frac{(t-2)!}{(a_1-1)!(a_2-1)!a_3!\cdots a_{10}!}+
\frac{(t-2)!}{(a_1-1)!a_2!(a_3-1)!a_4!\cdots a_{10}!}+\cdots \\
&+\frac{(t-2)!}{(a_1-1)!a_2!\cdots a_8!(a_9-1)!a_{10}!}+
\frac{(t-2)!}{(a_1-1)!a_2!\cdots a_9!(a_{10}-1)!}.
\end{align*}
This sum multiplied by $\frac{a_1!a_2!\cdots a_{10}!}{t!}$ 
simplifies to $\frac{a_1(t-a_1)}{t(t-1)}$, so
$\Pr[1]=a_1/(t-1)$. The same method works for $\Pr[2],\Pr[3],\ldots,\Pr[10]$.

\section*{Appendix B: Implementation}
We use the following C++ class to describe a deck of cards,
{\small
\begin{verbatim}
class Deck {
    Deck(...);
    double cardProb(int i); 
    double cardProb(int i, int d);
    void   removeCard(int i);
    void   addCard(int i); 
    ...
};
\end{verbatim} 
}

\noindent where $i=1,2,\ldots,10$ and $d=1,2,\ldots,10$.
The function  {\tt cardProb(i)} gives probabilities according to measure $\Prob$;
the function  {\tt cardProb(i,d)} gives probabilities according to measure $\Pr$,
when dealer's face-up card is $d$. The other functions are self explanatory.
The example below shows how to compute the probability of $\{1+2,\ 2+1\}$ for
cards coming from $\deck$.
{\small
\begin{verbatim}
double probability = 0.0;
double tmp = deck.cardProb(1);
if(tmp>0.0) {
    deck.removeCard(1);
    proability = 2.0*tmp*deck.cardProb(2);
    deck.addCard(1);
}
return(probability);
\end{verbatim}
}

Also it is convenient to create a class adding cards 
according to blackjack's rules. The class has an integer (hand's total) 
and a boolean variable (soft/hard hand) as members.

{\small
\begin{verbatim}
class Hand {
   Hand(...);
   Hand& operator+=(int i); 
   bool  operator<=(int i);
   bool  operator==(int i);
   bool  operator<(int i);
   ...
};

Hand operator+(const Hand& h, int i);
\end{verbatim}
}

The non-member operator {\tt  operator+(const Hand\& h, int i)} 
is used to create temporary objects passed by value to other functions;
given a hand, the operator  {\tt  operator+=(int i)} 
is used to add cards to it. For example, 
{\small
\begin{verbatim}
Hand h(5); h += 1;
\end{verbatim}
}

\noindent creates a soft sixteen.

With these two classes, 
the implementation of
the expected value functions  for stand (no naturals), double down, hit, and split
is straightforward.
{\small
\begin{verbatim}
double STAND(Hand p, int d, Deck deck) {
    if(p>21) return(-1.0);
    if(p<17) return(Q(23,d,deck)-Q(17,d,deck)-...-Q(21,d,deck));
    else return(1.0-Q(p,d,deck)-2.0*(Q(p+1,d,deck)+...+Q(21,d,deck)));       
}
\end{verbatim}
}
\noindent where {\tt Q(i,d,deck)} denotes $\Pr[D_d=i|\deck]$.

{\small
\begin{verbatim}
double DOUBLE(Hand p, int d, Deck deck) {
   double total = 0.0;
   for(int i=1;i<=10;++i) {
      double tmp = deck.cardProb(i,d);
      if(tmp>0.0) {
         deck.removeCard(i);
         total += tmp*STAND(p+i,d,deck); 
         deck.addCard(i);
      }  
   }
   return(2.0*total);
}
\end{verbatim}
}

A variable {\tt rec} controls the depth of recursion.
{\small
\begin{verbatim}
double HIT(Hand p, int d, Deck deck, int rec=13) {
   double total = 0.0;
   for(int i=1;i<=10;++i) {
      double tmp = deck.cardProb(i,d);
      if(tmp>0.0) {
         deck.removeCard(i);
         if(p+i>=21 || rec<=0) total += tmp*STAND(p+i,d,deck); 
         else total += tmp*max(STAND(p+i,d,deck),HIT(p+i,d,deck,--rec));
         deck.addCard(i);
      } 
   }
   return(total);
}
\end{verbatim}
}

When splitting we distinguish between a pair of aces and a non-ace pair.
With two aces we use  a boolean variable rsa (re-split aces).
{\small
\begin{verbatim}
RealNum SPLIT_ACES(int d, Deck deck, bool rsa) {
   Hand ace(1);
   double total = 0.0;
   for(int i=1;i<=10;++i) {  
      double tmp = deck.cardProb(i,d);
      if(tmp>0.0) {
         deck.removeCard(i);
         if(rsa && i==1) total += tmp*max(STAND(ace+i,d,deck),
                                          SPLIT_ACES(d,deck,false));
         else total += tmp*STAND(ace+i,d,deck);
         deck.addCard(i);
      }
   }
   return(2.0*total);
}
\end{verbatim}
}

Given $p>1$, we use two boolean variables, 
das (double down after split) and rsp (re-split non-ace pairs), and
an auxiliary procedure {\tt NOSPLIT(Hand p, int d, Deck deck, bool das).}
With {\tt rsp=true}, one hand could lead to four hands.
{\small 
\begin{verbatim}
double SPLIT(Hand p, int d, Deck deck, bool das, bool rsp) {
   double total = 0.0;
   for(int i=1;i<=10;++i) {
      double  tmp = deck.cardProb(i,d);
      if(tmp>0.0) {
         deck.removeCard(i);
         if(rsp && p==i) total += tmp*max(NOSPLIT(p+i,d,deck,das),
                                          SPLIT(p,d,deck,das,false));
         else total += tmp*NOSPLIT(p+i,d,deck,das);
         deck.addCard(i);
      }
   }
   return(2.0*total);
}

double NOSPLIT(Hand p, int d, Deck deck, bool das) {
   if(p>=21) return(STAND(p,d,deck));
   else
      if(das) return(max(STAND(p,d,deck),DOUBLE(p,d,deck),HIT(p,d,deck)));
      else return(max(STAND(p,d,deck),HIT(p,d,deck)));
}
\end{verbatim}
}

\section*{Appendix C: dealer's probabilities}
We start with a procedure that
computes dealer's probabilities of $17,18,\ldots,21$.
To be more specific let's focus on $\Prob[D_d=18]$; $19$, $20$, and $21$
work the same way, with 17 as an exception.

\subsubsection*{Computing: $\Prob[D_d=e]$}
We decompose $\{D_c=18\}$ into disjoint
subsets $\{d\to 16+2=18\}$, $\{d\to 15+3=18\}$,
$\{d\to 14+4=18\}$, and so on, where the notation $d\to e+k$ 
means all the possible cards combinations that add up to $e$ starting with
total $d$ and  exactly one card with face value $k$. So
$$
\Prob[D_d=18]=\Prob[d\to 16+2=18]+\Prob[d\to 15+3=18]+\cdots.
$$
Moreover, we notice that we can move backward, namely
$$
\Prob[d\to 16+2=18]=\Prob[2]\Prob[d\to 16|\deck-\{2\}].
$$
The implementation looks as follows.
{\small
\begin{verbatim}
double dealerProb21(int d, int e, Deck deck) {
   double total = 0.0;
   if(e==17)  return(dealerProb17(d,17,deck));
   else {
     for(int i=e-16;i<=min(e-d,11);++i) {
        double tmp = deck.cardProb(i);
        if(tmp>0.0) {
           deck.removeCard(i);
           total += tmp*dealerProb17(d,e-i,deck);
           deck.addCard(i);
        }
     }
   }
   return(total);
}
\end{verbatim}
}
 
Now, a procedure {\tt dealerProb17(d,e)} that computes 
the probabilities $\Prob[d\to e]$ for $e\le 17$ is needed. 

\subsubsection*{Computing: $\Prob[d\to e\mbox{ soft}]$ and $\Prob[d\to e\mbox{ hard}]$}
Given $e\le 17$,
we distinguish between soft and hard totals:
$\Prob[d\to e]= \Prob[d\to e\mbox{ soft}]+\Prob[d\to e\mbox{ hard}]$, or
\begin{verbatim}
   dealerProb17(d,e) = dealerSoft17(d,e) + dealerHard17(d,e);
\end{verbatim}
A brute force is used to compute $\Prob[d\to e\mbox{ soft}]$;
the most complicated case $2\to 17s$ can be decomposed into five
cases: $2+1+(4)$, $2+1+1+(3)$, $2+(2)+1+(2)$, $2+(3)+1+1$, and $2+(4)+1$,
where $(a)$ means all the decompositions of $a$; the procedure
{\tt prob(int a)} computes the probability of them.

Next, $\Prob[d\to e\mbox{ hard}]$ is computed using a recursive procedure
presented below.
{\small
\begin{verbatim}
double dealerHard17(int d, int e, Deck deck) {
   if(e==d) return(1.0);
   if(e-d==1) return(0.0);
   double total = 0.0;
   double tmp = deck.cardProb(1);
   if(tmp>0.0) {
      deck.removeCard(1);
      if(d+11<=16) total += tmp*dealerSoft2Hard17(d+11,e,deck);
      deck.addCard(1);
   }
   for(int i=2;i<=min(e-d,10);++i) {
      double tmp = deck.cardProb(i);
      if(tmp>0.0) {
         deck.removeCard(i);
         if(d+i<11) total += tmp*dealerHard17(d+i,e,deck);
         else total += tmp*prob(e-d-i,deck);
         deck.addCard(i);
      }
   }
   return(total);
}
\end{verbatim}
}

\noindent And we use a brute force to compute {\tt dealerSoft2Hard17(d,e,deck)},
for example, $16s\to 15h$ can be decomposed as 
$16s+9$, $16s+8+1$, $16s+7+(2)$, and $16s+6+(3)$.

Finally, two auxiliary procedures are needed. The first one, {\tt prob(int a)}, computes
the probability of all possible decompositions of $a$. For example,
if $a=4$ then there are eight of them: $4$, $1+3$, $3+1$, $2+2$, $1+1+2$,
$1+2+1$, $2+1+1$, and $1+1+1+1$. 
Since we need it only for $a=0,1,2,\ldots,6$, 
we can list all the cases and compute them one by one.

The second procedure, {\tt prob(int a, int b)}, where $a,b\le 5$, 
computes the probability of all possible 
decompositions of $a$ and $b$. For example, if $a=2$ and $b=3$, then
$2+3$, $1+1+3$, $2+1+2$, $2+2+1$, $1+1+1+2$, $1+1+2+1$, $2+1+1+1$,
and $1+1+1+1+1$. 
{\small
\begin{verbatim}
double prob(int a, int b, Deck deck) {
   if(a==0) return(prob(b,deck));
   if(b==0) return(prob(a,deck));
   double total = 0.0;
   for(int i=1;i<=a;++i) {
      double tmp=deck.cardProb(i);
      if(tmp>0.0) {
         deck.removeCard(i);
         total += tmp*prob(a-i,b,deck);
         deck.addCard(i);
      }
   }
   return(total);
}
\end{verbatim}
}

\end{document}